\documentstyle[12pt, page1]{article}
\newcommand{\sect}[1]{\section{#1}\setcounter{equation}{0}}

\font\mbn=msbm10 scaled \magstep1
\font\mbs=msbm7 scaled \magstep1
\font\mbss=msbm5 scaled \magstep1
\newfam\mbff
\textfont\mbff=\mbn
\scriptfont\mbff=\mbs
\scriptscriptfont\mbff=\mbss
\def\mbf{\fam\mbff}
\def\Re{{\mbf R}}

\def\Z{{\mbf Z}}
\def\Co{{\mbf C}}

\def\Di{{\mbf D}}
\newtheorem{Th}{Theorem}[section]
\newtheorem{Lm}[Th]{Lemma}

\newtheorem{Proposition}[Th]{Proposition}

\author{Alexander Brudnyi\\
Department of Mathematics and Statistics\\
University of Calgary\\
Calgary, Canada}
\title{A Note on the Geometry of Green-Lazarsfeld Sets}
\date{}
\begin{document}
\maketitle
\begin{abstract}
{We study some subsets of the Green-Lazarsfeld set $\Sigma^{1}(M)$ for a 
compact K\"{a}hler manifold $M$.}
\end{abstract}
\sect{\hspace*{-1em}. Introduction.} 
\noindent  
Let $M$ be a compact K\"{a}hler manifold and 
$H^{1}(M,\Co^{*}):=Hom(\pi_{1}(M),\Co^{*})$. Given a character 
$\xi\in H^{1}(M,\Co^{*})$ let $\Co_{\xi}$ denote the associated 
$\pi_{1}(M)$-module. Let $\Sigma^{1}(M)$ be the set of characters
$\xi$ such that $H^{1}(\pi_{1}(M),\Co_{\xi})$ is nonzero.
The structure of $\Sigma^{1}(M)$ was described in the papers
of Beauville [Be], Simpson [S], and Campana [Ca]:

{\em There is a finite number of surjective
holomorphic maps with connected fibres $f_{i}:M\rightarrow C_{i}$
onto smooth compact complex curves of genus $g\geq 1$
and torsion characters $\rho_{i},\xi_{j}
\in H^{1}(M,\Co^{*})$ such that}
$$
\Sigma^{1}(M)=\bigcup_{i}\rho_{i}f_{i}^{*}H^{1}(C_{i},\Co^{*})\cup
\bigcup_{j}\{\xi_{j}\}\ .
$$
In this paper we study the components $f_{i}^{*}H^{1}(C_{i},\Co^{*})$
corresponding to trivial $\rho_{i}$. 
\begin{Th}\label{te1}
Let $f: M\rightarrow C$ be a surjective holomorphic map with connected 
fibres of a compact K\"{a}hler manifold $M$ onto a smooth 
complex curve $C$ of genus $g\geq 1$.
Assume that $\xi\in H^{1}(C,\Co^{*})$ is not torsion. Then 
$$
f^{*}:H^{1}(\pi_{1}(C),\Co_{\xi})\longrightarrow 
H^{1}(\pi_{1}(M),\Co_{f^{*}\xi})
$$
is an isomorphism.
\end{Th}
\sect{\hspace*{-1em}. Proof of Theorem 1.1.}
\noindent {\bf 2.1.} 
Let $T_{2}(\Co)\subset GL_{2}(\Co)$ be the Lie subgroup of
upper triangular matrices. By $T_{2}\subset T_{2}(\Co)$ we denote
the Lie group of matrices of the form
\[
\left(
\begin{array}{cc}
a&b\\ 0&1\\
\end{array}
\right)\ \ \ \ \ (a\in\Co^{*},\ b\in\Co) .
\]
Also by $D_{2}\subset T_{2}$ we denote the subgroup of diagonal matrices
and by $N_{2}\subset T_{2}$ the subgroup of unipotent matrices.
Let $M$ be a compact K\"{a}hler manifold.
For a homomorphism $\rho\in Hom(\pi_{1}(M),T_{2})$ we let
$\rho_{a}\in Hom(\pi_{1}(M),\Co^{*})$ denote the upper diagonal character
of $\rho$.
Further, $N_{2}$ acts on $T_{2}$ by conjugation. Any two
homomorphisms from $Hom(\pi_{1}(M),T_{2})$ belonging to the orbit
of this action will be called {\em equivalent}.
It is well known that the class of equivalence of 
$\rho$ is uniquely defined by an element 
$c_{\rho}\in H^{1}(\pi_{1}(M),\Co_{\rho_{a}})$ (see e.g. [A, Prop. 2]).
Now our theorem is a consequence of the following result.
\begin{Th}\label{curv}
Assume that $\rho\in Hom(\pi_{1}(M),T_{2})$ 
satisfies $Ker(f_{*})\subset Ker(\rho_{a})$ but $Ker(f_{*})\not\subset 
Ker(\rho)$. Then $\rho_{a}$ is a torsion character.
\end{Th}
The next several sections contain some results used in the proof of 
Theorem \ref{curv}.\\
{\bf 2.2.} Let $f:M\rightarrow C$ be a holomorphic surjective map with
connected fibres of a compact K\"{a}hler manifold $M$ onto a smooth compact
complex curve $C$. Consider a flat vector bundle $L$ on $C$ of complex rank 
1 with unitary structure group and let $E=f^{*}L$ be the pullback of this 
bundle to $M$. Let $\omega\in\Omega^{1}(E)$
be a holomorphic 1-form with values in $E$ 
(by the Hodge decomposition $\omega$ is $d$-closed).
\begin{Lm}\label{vanish}
Assume that $\omega|_{V_{z}}=0$ for the fibre $V_{z}:=f^{-1}(z)$ over a 
regular value $z\in C$. Then $\omega|_{V}=0$ for any fibre $V$ of $f$.
\end{Lm}
If $V$ is not smooth the lemma asserts that $\omega$ equals 0 on the 
smooth part of $V$.\\
{\bf Proof.} Denote by $S\subset C$ the finite set
of non-regular values of $f$. By Sard's theorem, 
$f:M\setminus f^{-1}(S)\rightarrow C\setminus S$ is a fibre bundle with 
connected fibres. According to our assumption there is a 
fibre $V_{z}$ of this bundle such that $\omega|_{V_{z}}=0$. 
Then $\omega|_{V}=0$ for any fibre $V$ of $f|_{M\setminus f^{-1}(S)}$.
In fact, for any fibre $V$ there is an open neighbourhood $O_{V}$ of $V$
diffeomorphic to $\Re^{2}\times V$ such
that $E|_{O_{V}} (\cong E|_{V})$ is the trivial flat vector bundle (here
$E|_{V}$ is trivial as the pullback of a bundle defined over a 
point). Any $d$-closed holomorphic 1-form defined on $O_{V}$ that vanishes on
$V$ is $d$-exact, and so it vanishes on each fibre contained in $O_{V}$. 
Starting with a tubular neighbourhood $O_{V_{z}}$ and taking into account 
that $\omega|_{O_{V_{z}}}$ can be considered as a $d$-closed holomorphic 
1-form (because $E|_{O_{V_{z}}}$ is trivial) we obtain
that $\omega|_{V}=0$ for any fibre $V\subset O_{V_{z}}$. Then the required
statement follows by induction if we cover $M\setminus f^{-1}(S)$ by open 
tubular neighbourhoods of fibres of $f$ and use the fact that 
$C\setminus S$ is connected.

Consider now fibres over the singular part $S$.  Let $O_{x}\subset C$ be the
neighbourhood of a point $x\in S$ such that 
$O_{x}\setminus S=O_{x}\setminus\{x\}$ is biholomorphic to 
$\Di\setminus\{0\}$. In particular, $\pi_{1}(O_{x}\setminus S)=\Z$. 
Without loss of generality we may assume that $V=f^{-1}(x)$ is a deformation 
retract of an open neighbourhood $U_{x}$ of $V$ and 
$W_{x}:=f^{-1}(O_{x})\subset U_{x}$ (such $U_{x}$ exists, e.g.,
by the triangulation theorem of the pair $(M,V)$, see [L, Th.2]). Since 
$E|_{V}$ is trivial, the bundle $E|_{U_{x}}$ is also trivial. Therefore we 
can regard $\omega$ as a $d$-closed holomorphic 1-form on $U_{x}$. Moreover, 
we already have proved that $\omega|_{F}=0$ for any fibre 
$F\subset W_{x}\setminus f^{-1}(S)$. So $\omega$ equals 0 on each fibre of 
the fibration $f:W_{x}\setminus f^{-1}(S)\rightarrow O_{x}\setminus S$. 
This implies that there is a $d$-closed
holomorphic  1-form $\omega_{1}$ on $O_{x}\setminus S$ such that $\omega=
f^{*}(\omega_{1})$. Assume now that $\omega|_{V}\neq 0$. Then integration
of $\omega|_{W_{x}}$ along paths determines a non-trivial homomorphism 
$h_{\omega,x}:\pi_{1}(W_{x})\rightarrow\Co$ whose image
is isomorphic to $\Z$. Indeed, 
since embedding $W_{x}\setminus f^{-1}(S)\hookrightarrow W_{x}$ induces a 
surjective homomorphism of fundamental groups, we can integrate $\omega$ by 
paths contained in $W_{x}\setminus f^{-1}(S)$ obtaining the same image in 
$\Co$. Further, for any path $\gamma\subset W_{x}\setminus f^{-1}(S)$ we have
$\int_{\gamma}\omega=\int_{\gamma_{1}}\omega_{1}$, where $\gamma_{1}$ is a 
path in $O_{x}\setminus S$ representing the element 
$f_{*}(\gamma)\in\pi_{1}(O_{x}\setminus S)$. But 
$\pi_{1}(O_{x}\setminus S)=\Z$ and therefore 
$h_{\omega,x}(\pi_{1}(W_{x}))\cong\Z$. Now, any path in $W_{x}$ is 
homotopically equivalent inside $U_{x}$ to a path contained in $V$ and so 
we can define $h_{\omega,x}$ integrating $\omega|_{V}$. Thus we
obtain a homomorphism $h$ of 
$\Gamma:=(\pi_{1}(V)/D\pi_{1}(V))/_{torsion}$ into $\Co$ with the same
image as for $h_{\omega,x}$, i.e., $\cong\Z$. 
Without loss of generality we may assume
that $V$ is smooth (for otherwise, we apply the arguments below to each
irreducible component of a desingularization of $V$). Integrating holomorphic
1-forms on $V$ along paths we embed  $\Gamma$ into some $\Co^{k}$ as a 
lattice of rank $2k$. Then there is a linear holomorphic functional 
$f_{\omega}$ on $\Co^{k}$ such that $h=f_{\omega}|_{\Gamma}$. Since 
$rk(h(\Gamma))=1$, $f_{\omega}$ equals 0 on a subgroup $H\subset\Gamma$ 
isomorphic to $\Z^{2k-1}$. In particular,
$f_{\omega}=0$ on the vector space $H\otimes\Re$ of real dimension $2k-1$. 
But $Ker(f_{\omega})$ is a complex vector space. So $f_{\omega}=0$ on 
$\Co^{k}$ 
which implies that $\omega|_{V}=0$. This contradicts our assumption and 
proves the required statement for fibres over the points of $S$.

The lemma is proved.\ \ \ \ \ $\Box$
\\
{\bf 2.3.} Next we prove that the character $\rho_{a}$ in 
Theorem \ref{curv} cannot be non-unitary.
\begin{Proposition}\label{normsubgr}
Let $\rho\in Hom(\pi_{1}(M),T_{2})$ be such that $\rho_{a}$ is a non-unitary
character. Assume that $\rho_{a}|_{Ker(f_{*})}$ is trivial. Then
$\rho|_{Ker(f_{*})}$ is also trivial.
\end{Proposition}
{\bf Proof.}
Let $p_{t}:M_{t}\rightarrow M$ be the Galois covering with  
transformation group $Tor(\pi_{1}(M)/D\pi_{1}(M))$. Then 
$\rho|_{\pi_{1}(M_{t})}$ determines a $C^{\infty}$-trivial complex rank-2 
flat vector bundle on $M_{t}$, because 
$\rho_{a}|_{\pi_{1}(M_{t})}=\exp(\rho'|_{\pi_{1}(M_{t})})$ for some
$\rho'\in Hom(\pi_{1}(M),\Co)$. In particular, $\rho|_{\pi_{1}(M_{t})}$
can be defined by a flat connection 
\[
\Omega=\left(
\begin{array}{cc}
\omega&\eta\\0&0\\
\end{array}
\right);\ \ \ \ \ \ \  d\Omega-\Omega\wedge\Omega=0
\]
on $M_{t}\times\Co^{2}$ where $\omega$ is a $d$-harmonic 1-form on $M_{t}$
lifted from $M$.
Moreover, since $\rho_{a}|_{Ker(f_{*})}$ is trivial, $\omega$ is the
pullback by $f\circ p_{t}$ of a $d$-harmonic 1-form defined on $C$. Let
$\omega=\omega_{1}+\omega_{2}$ be the type decomposition of $\omega$ into the
sum of holomorphic and antiholomorphic 1-forms lifted from $C$. Denote by 
$E_{\rho}$ the rank-1 flat vector bundle on $M_{t}$ with unitary structure 
group constructed by the flat connection $\omega_{2}-\overline{\omega_{2}}$ 
and by $E_{0}=M_{t}\times\Co$ the trivial flat vector bundle. Observe that 
$E_{\rho}$ is pullback of a flat bundle on $C$. In particular, 
$E_{\rho}|_{V}$ is trivial for any fibre $V$ of $f\circ p_{t}$.  According 
to [Br, Th.1.2] the equivalence
class of $\rho|_{\pi_{1}(M_{t})}$ is determined by a $d$-harmonic 1-form 
$\theta$
with values in $End(E_{\rho}\oplus E_{0})$ satisfying $\theta\wedge\theta$
represents 0 in $H^{2}(M_{t},End(E_{\rho}\oplus E_{0}))$. More precisely,
\[
\theta=
\left(
\begin{array}{cc}
\omega_{1}&\eta'\\0&0\\
\end{array}
\right)
\]
where $\eta'$ is a $d$-harmonic 1-form with values in $E_{\rho}$ satisfying
$\omega_{1}\wedge\eta'$ represents 0 in $H^{2}(M_{t},E_{\rho})$.
Let $M_{t}\stackrel{g_{1}}{\longrightarrow}M_{1}\stackrel{g_{2}}
{\longrightarrow}C$ be the Stein factorization of 
$f\circ p_{t}$. Here $g_{1}$ is a morphism with connected
fibres onto a smooth curve $M_{1}$ and $g_{2}$ is a finite morphism. 
Then for any fibre $V\hookrightarrow
M_{t}$ of $g_{1}$ we have that $p_{t}|_{V}:V\rightarrow p_{t}(V)$ is a
regular covering of the fibre $p_{t}(V)$ of $f$ with a finite abelian
transformation group.
\begin{Lm}\label{van}
$\eta'|_{V}=0$ for fibres over regular values of $g_{1}$.
\end{Lm}
Based on this statement we, first, finish the proof of the proposition and 
then will prove the lemma. From the lemma and Lemma \ref{vanish} it follows 
that $\eta'|_{W}=0$ for any fibre $i:W\hookrightarrow M_{t}$ of $g_{1}$.  Thus
$\rho|_{\pi_{1}(M_{t})}$ is trivial on $i_{*}(\pi_{1}(W))\subset
\pi_{1}(M_{t})$. But $\pi_{1}(W)$ is a subgroup of a finite index in
$\pi_{1}(p_{t}(W))$ and the image of $\rho|_{j_{*}(\pi_{1}(p_{t}(W)))}$
consists of unipotent matrices by the assumption of the proposition. (Here
$j:p_{t}(W)\hookrightarrow M$ is embedding.) Therefore 
$\rho|_{j_{*}(\pi_{1}(p_{t}(W)))}$ is trivial for any $W$.
Denote by $E$ the flat vector bundle on $M$ associated to $\rho$. Then 
we have proved that $E|_{W}$ is trivial for any fibre $W$ of $f$. Further, 
let $(U_{i})_{i\in I}$ be an open
cover of $C$ such that $W_{i}:=f^{-1}(U_{i})$ is an open neighbourhood of a
fibre $V_{i}=f^{-1}(x_{i})$, $x_{i}\in U_{i}$, and $W_{i}$ is deformable onto
$V_{i}$. Since $E|_{V_{i}}$ is trivial, $E|_{W_{i}}$ is also trivial. In
particular, $E$ is defined by a locally constant cocycle $\{c_{ij}\}$ 
defined on the cover $(W_{i})_{i\in I}$. But then $\{c_{ij}\}$ is the 
pullback of a cocycle defined on $(U_{i})_{i\in I}$ because the fibres of
$f$ are connected. This cocycle determines a bundle $E'$ on $C$ such that 
$f^{*}E'=E$. Then $\rho$ is the pullback of a homomorphism
$\rho'\in Hom(\pi_{1}(C),T_{2})$ constructed by $E'$.

This completes the proof of the proposition modulo Lemma \ref{van}.\\ 
{\bf Proof of Lemma \ref{van}.} 
Observe that $\omega_{1}\neq 0$ in the above definition of $\theta$ 
because $\rho_{a}$ is non-unitary. We also regard $\eta'\neq 0$, 
for otherwise, the image of $\rho$ consists of diagonal matrices and the 
required statement is trivial. Let
$i:V\hookrightarrow M_{t}$ be a fibre over a regular value of $g_{1}$. For
any $\lambda\in\Co^{*}$ consider the form
\[
\theta_{\lambda}=\left(
\begin{array}{cc}
\lambda\omega_{1}&\eta'\\0&0\\
\end{array}
\right).
\]
Clearly, $\theta_{\lambda}\wedge\theta_{\lambda}$ represents 0 in
$H^{2}(M_{t},End(E_{\rho}\oplus E_{0}))$ and thus, by 
[Br, Prop.2.4], it determines a
representation $\rho_{\lambda}\in Hom(\pi_{1}(M_{t}),T_{2})$ with the upper
diagonal character $\rho_{\lambda a}$ defined by the flat connection
$\lambda\omega_{1}+\omega_{2}$. In particular, the family
$\{\rho_{\lambda a}\}$ contains infinitely many different characters. 
Assume that $\rho|_{i_{*}(\pi_{1}(V))}$ is not trivial. Then it can be 
determined by the restriction  $\theta|_{i(V)}$. But according to our 
assumption
$\omega_{1}|_{i(V)}=0$ and $E_{\rho}|_{i(V)}$ is a trivial flat vector 
bundle. Thus $\rho|_{i_{*}(\pi_{1}(V))}$ is
defined by $\eta'|_{i(V)}\in H^{1}(i(V),\Co)$. Let $\psi\in
Hom(i(V),\Co)$ be a homomorphism obtained by integration of
$\eta'|_{i(V)}$ along paths generating $\pi_{1}(i(V))$. Then
\[
\rho|_{i_{*}(\pi_{1}(V))}=\left(
\begin{array}{cc}
1&\psi\\0&1\\
\end{array}
\right).
\]
We also obtain that $\rho_{\lambda}|_{i_{*}(\pi_{1}(V))}=\rho|_{
i_{*}(\pi_{1}(V))}$ and so it is not trivial for any $\lambda\in\Co^{*}$. 
We will prove now that the family $\{\rho_{\lambda a}\}$ contains finitely 
many different characters. This contradicts our assumption and shows that
$\rho|_{i_{*}(\pi_{1}(V))}$ is trivial and so $\eta'|_{V}=0$.

Let $H:=i_{*}(\pi_{1}(V))$ be a normal subgroup of $\pi_{1}(M_{t})$ (the 
normality follows from the exact homotopy sequence for the bundle defined
over regular values of $g_{1}$).
Denote by $H^{ab}$ the quotient $H/DH$. Then $H^{ab}$ is an 
abelian group of a finite rank. Moreover, $H^{ab}$ is a normal subgroup of 
$K:=\pi_{1}(M_{t})/DH$ and the group $K_{1}:=K/H^{ab}$ is finitely generated. 
Note that $\rho_{\lambda}$ induces a homomorphism 
$\widehat\rho_{\lambda}:K\rightarrow T_{2}$ for any
$\lambda\in\Co^{*}$ with the same image as for $\rho_{\lambda}$
because the image of $\rho_{\lambda}|_{H}$ is abelian. In particular,
$\widehat\rho_{\lambda}(H^{ab})$ is a normal subgroup of $\widehat\rho(K)$. 
Since by our assumption $\widehat\rho_{\lambda}(H^{ab})$ is a non-trivial 
subgroup of unipotent matrices, from the identity
\begin{equation}\label{ident}
\left(
\begin{array}{cc}
c&b\\0&1\\
\end{array}
\right)\cdot
\left(
\begin{array}{cc}
1&v\\0&1\\
\end{array}
\right)\cdot
\left(
\begin{array}{cc}
c&b\\0&1\\
\end{array}
\right)^{-1}=
\left(
\begin{array}{cc}
1&c\cdot v\\0&1\\
\end{array}
\right)\ \ \ \ \ (c\in\Co^{*},\ b,\ v\in\Co)
\end{equation}
it follows that the action of $\widehat\rho_{\lambda}(K)$ on
$\widehat\rho_{\lambda}(H^{ab})$ by
conjugation is defined by multiplication of non-diagonal elements of
$\widehat\rho_{\lambda}(H^{ab})$ by elements of 
$\rho_{\lambda a}(\pi_{1}(M_{t}))$. Note also that $Tor(H^{ab})$ belongs to
$Ker(\widehat\rho_{\lambda})$ for any $\lambda$.
Thus if $H'\cong\Z^{s}$ is a maximal free 
abelian subgroup of $H^{ab}$ then $\rho_{\lambda a}(\pi_{1}(M_{t}))$ 
consists of 
eigen values of matrices from $SL_{s}(\Z)$ obtained by the natural action 
(by conjugation) of $K_{1}$ on $H'$. Since $K_{1}$ is
finitely generated, the number of different $\rho_{\lambda a}$ is finite
which is false.

This contradiction completes the proof of the lemma.\ \ \  \ \ $\Box$
\\
{\bf 2.4.} Let $f:M\rightarrow C$ be a holomorphic map with connected 
fibres of a compact K\"{a}hler manifold $M$ onto a smooth
compact complex curve $C$ of genus $g\geq 1$. Then $\pi_{1}(M)$ is defined 
by the exact sequence
$$
\{1\}\longrightarrow Ker(f_{*})\longrightarrow\pi_{1}(M)\longrightarrow
\pi_{1}(C)\longrightarrow\{1\}\ . $$
Let $G\subset Ker(f_{*})$ be a normal subgroup of $\pi_{1}(M)$. The 
quotient $R:=\pi_{1}(M)/G$ is
defined by the sequence 
$$ 
\{1\}\longrightarrow Ker(f_{*})/G\longrightarrow
R\longrightarrow\pi_{1}(C)\longrightarrow\{1\}\ . 
$$ 
Assume that

(1)\ \ $H:=Ker(f_{*})/G$ is a free abelian group of finite rank $k\geq 1$;

(2)\ \ the natural action of $D\pi_{1}(C)$ on $H$ is trivial.
\\
From (2) it follows that the action $s$ of $\pi_{1}(C)$ on $H$ determines an 
action of $\Z^{2g}\cong\pi_{1}(C)/D\pi_{1}(C)$ on $H$. Identifying $H$
with $\Z^{k}$ we can think of $H$ as a subgroup (lattice) of 
$Z^{k}\otimes\Co=\Co^{k}$. Then $s$ determines a representation 
$s':\pi_{1}(C)\rightarrow GL_{k}(\Co)$ such that 
$D\pi_{1}(C)\subset Ker(s')$ and 
$s'(g)|_{H}=s(g)$ for any $g\in\pi_{1}(C)$. 
Since $s'$ descends to a representation $\Z^{2g}\rightarrow GL_{k}(\Co)$,
it admits a decomposition $s'=\oplus_{j=1}^{m}s_{j}$ where $s_{j}$ is
equivalent to a nilpotent representation 
$\pi_{1}(C)\rightarrow T_{k_{i}}(\Co)$ with a 
diagonal character $\rho_{j}$. Here $\sum_{j=1}^{m}k_{j}=k$.
\begin{Proposition}\label{curve}
All characters $\rho_{j}$ are torsion.
\end{Proposition}
{\bf Proof.} By definition, $R$ is an extension of $\pi_{1}(C)$ by $H$. It 
is well known, see, e.g. [G, Ch.I, Sec.6], that the class of extensions 
equivalent to $R$ is uniquely defined by an element 
$c\in H^{2}(\pi_{1}(C),H)$ where the cohomology is defined by the action
$s$ of $\pi_{1}(C)$ on $H$. Let $f\in Z^{2}(\pi_{1}(C),H)$ be a cocycle 
determining $c$. Then one can define a representative of the 
equivalence class of extensions as the direct product $H\times\pi_{1}(C)$ 
with 
multiplication 
$$
\begin{array}{c}
(h_{1},g_{1})\cdot (h_{2},g_{2})=
(h_{1}+s(g_{1})(h_{2})+f(g_{1},g_{2}),g_{1}\cdot g_{2});\\ 
\\
h_{1},h_{2}\in H,\ g_{1},g_{2}\in\pi_{1}(C)\ .
\end{array}
$$ 
The natural embedding $H\hookrightarrow\Z^{k}\otimes\Co (=\Co^{k})$
determines an embedding $i$ of $R$ into the group $R'$ defined as
$\Co^{k}\times\pi_{1}(C)$ with multiplication 
$$
\begin{array}{r}
(v_{1},g_{1})\cdot (v_{2},g_{2})=
(v_{1}+s'(g_{1})(v_{2})+f(g_{1},g_{2}),g_{1}\cdot g_{2});\\
\\
f\in Z^{2}(\pi_{1}(C),H),\ v_{1},v_{2}\in\Co^{k},\ g_{1},g_{2}\in\pi_{1}(C)\ .
\end{array}
$$ 
Here we regard $f$ as an element of $Z^{2}(\pi_{1}(C),\Co^{k})$ defined by
the action $s'$. From the decomposition $s'=\oplus_{j=1}^{m}s_{j}$ it follows
that there is an invariant $\pi_{1}(C)$-submodule $V_{j}\subset\Co^{k}$ of
$dim_{\Co}V_{j}=k-1$ such that $W_{j}=\Co^{k}/V_{j}$ is a one-dimensional
$\pi_{1}(C)$-module and the action of $\pi_{1}(C)$ on $W_{j}$ is defined as
multiplication by the character $\rho_{j}$. Then, by definition, $V_{j}$ is 
a normal subgroup of $R'$ and the the quotient group $R_{j}=R'/V_{j}$ is 
defined by the sequence 
$$
\{1\}\longrightarrow\Co\longrightarrow
R_{j}\longrightarrow\pi_{1}(C)\longrightarrow\{1\}\  .
$$ 
Here the action of $\pi_{1}(C)$ on $\Co$ is multiplication by the character 
$\rho_{j}$. Further, the equivalence class of extensions isomorphic to 
$R_{j}$ is defined by an element $c_{j}\in H^{2}(\pi_{1}(C),\Co)$ (the 
cohomology is defined by the above action of $\pi_{1}(\Co)$ on $\Co$). We
will assume that the character $\rho_{j}$ is non-trivial (for otherwise,
$\rho_{j}$ is clearly torsion). Let us
denote by $t_{j}$ the composite homomorphism $\pi_{1}(M)\longrightarrow
R\stackrel{i}{\longrightarrow}R'\longrightarrow R_{j}$.

Let $E_{\rho_{j}}$ be a complex rank-1 flat vector bundle on $C$ constructed 
by $\rho_{j}\in Hom(\pi_{1}(C),\Co^{*})$. Since $C$ is a 
$K(\pi_{1}(C),1)$-space, there is a natural isomorphism of the above group 
$H^{2}(\pi_{1}(C),\Co)$ and the \v{C}ech cohomology group 
$H^{2}(C, {\bf E_{\rho_{j}}})$ of the sheaf of locally
constant sections of $E_{\rho_{j}}$, see e.g. [M, Ch.1, Complement to Sec.2].
But each flat vector bundle on $C$ is $C^{\infty}$-trivial and each 
homomorphism from $Hom(\pi_{1}(C),\Co^{*})$ can be continuously deformed 
inside of $Hom(\pi_{1}(C),\Co^{*})$ to the trivial homomorphism. Therefore 
by the index theorem 
$$ 
dim_{\Co}H^{0}(C,{\bf E_{\rho_{j}}})-dim_{\Co}H^{1}(C,{\bf E_{\rho_{j}}})+
dim_{\Co} H^{2}(C,{\bf E_{\rho_{j}}})=\chi(\pi_{1}(C))=2-2g\ .
$$
Note that $H^{0}(C,{\bf E_{\rho_{j}}})=0$ because $\rho_{j}$ is non-trivial.
Furthermore, $H^{1}(C,{\bf E_{\rho_{j}}})$ is in a one-to-one correspondence 
with the set of non-equivalent representations 
$\rho:\pi_{1}(C)\longrightarrow T_{2}\subset T_{2}(\Co)$ with the
upper diagonal character $\rho_{j}$. Using the identity
$\prod_{i=1}^{g}[e_{i},e_{g+i}]=e$ for generators 
$e_{1},..,e_{2g}\in\pi_{1}(C)$ we easily obtain that 
$dim_{\Co}H^{1}(C,{\bf E_{\rho_{j}}})=2g-2$. Thus we have
$H^{2}(C,{\bf E_{\rho_{j}}})=0$. This shows that $H^{2}(\pi_{1}(C),\Co)=0$ and
$R_{j}$ is isomorphic to the semidirect product of $\Co$ and $\pi_{1}(C)$, 
i.e., $R_{j}=\Co\times\pi_{1}(C)$ with multiplication 
$$
(v_{1},g_{1})\cdot (v_{2},g_{2})=
(v_{1}+\rho_{j}(g_{1})\cdot v_{2},g_{1}\cdot g_{2}),\ \ \
v_{1},v_{2}\in \Co,\ g_{1},g_{2}\in\pi_{1}(C)\ .
$$ 
Let us determine a map $\phi_{j}$ of $R_{j}$ to $T_{2}$
by the formula
\[
\phi_{j}(v,g)=\left(
\begin{array}{cc}
\rho_{j}(g)&v\\ 0&1\\
\end{array}
\right)
\]
Obviously, $\phi_{j}$ is a correctly defined homomorphism with  upper 
diagonal character $\rho_{j}$. Hence 
$\phi_{j}\circ t_{j}:\pi_{1}(M)\rightarrow T_{2}$ is a homomorphism 
which is non-trivial on $Ker(f_{*})$ by its definition. Now Proposition 
\ref{normsubgr} implies that $\rho_{j}$ is a unitary
character. Therefore we have proved that all characters of the action of 
$\pi_{1}(C)$ on $H$ are unitary. This means that each element of the above 
action is defined by a matrix from $SL_{k}(\Z)$ with unitary eigen values. 
Applying to these eigen values the Kronecker theorem asserting that an 
algebraic integer is a root of unity if and only if all its Galois conjugates
are unitary (see, e.g. [BS, p.105,Th.2]) we obtain that each $\rho_{j}$ is a 
torsion character.\ \ \ \ \ $\Box$\\
{\bf 2.5.} In this section we will prove a result
about the structure of $Ker(f_{*})$ where $f:M\rightarrow C$ is a 
surjective holomorphic map with connected fibres of a compact K\"{a}hler 
manifold $M$ onto a curve $C$.

Using the Lojasiewicz triangulation theorem [L, Th.2], compactness of 
$M$ and the fact that all fibres of $f$ are connected, we can find an
open finite cover $(U_{i})_{1\leq i\leq s}$ of $C$ such that 
$W_{i}:=f^{-1}(U_{i})$
is an open neighbourhood of a fibre $V_{i}=f^{-1}(x_{i})$, $x_{i}\in U_{i}$,
and $V_{i}$ is a deformation retract of some open 
$\widetilde W_{i}\supset W_{i}$. Let us fix some points $x_{i}\in V_{i}$
and paths $\gamma_{i}$ connecting $x_{1}$ with $x_{i}$, $1\leq i\leq s$.
By $G_{i}$ we denote the image of 
$\gamma_{i}\cdot\pi_{1}(V_{i},x_{i})\cdot\gamma_{i}^{-1}$ in 
$\pi_{1}(M,x_{1})$. (Clearly, each $G_{i}$ is finitely generated.)
\begin{Lm}\label{kern}
$Ker(f_{*})$ is the minimal normal subgroup of $\pi_{1}(M,x_{1})$
containing all $G_{i}$ $(1\leq i\leq s)$.
\end{Lm}
{\bf Proof.} Let $R\subset\pi_{1}(M,x_{1})$ be the minimal normal subgroup
containing all $G_{i}$. Clearly, $R\subset Ker(f_{*})$. We will prove the 
reverse inclusion. Let $p:M_{1}\rightarrow M$ be the regular covering over 
$M$ with the transformation group $R_{1}:=\pi_{1}(M)/R$. It is a principle 
bundle over $M$ with discrete fibre $R_{1}$. We will prove that the 
restriction of $M_{1}$ to each $\widetilde W_{i}$ is the trivial bundle. 
Since $V_{i}$ is a deformation retract of $\widetilde W_{i}$,
it suffices to prove that for any closed path $\gamma\subset V_{i}$ 
passing through $x_{i}$ and any point $y\in p^{-1}(x_{i})$ the unique 
covering path $\gamma'\subset M_{1}$ of $\gamma$ passing through $y$ is 
closed. To prove this let us consider the path
$s:=\gamma_{i}\cdot\gamma\cdot\gamma_{i}^{-1}$. Let 
$s_{1}\subset M_{1}$ be any path which covers $s$. Then $s_{1}$ is closed,
because $s$ represents an element of the {\em normal subgroup }
$R\subset\pi_{1}(M)$. Let $\gamma'$ be the covering
of $\gamma$ with endpoints $y$ and $y'$. Then there is an element 
$g\in R_{1}$ such that $y'=g(y)$. Let $h$ and $h'$ be the unique
covering paths of $\gamma_{i}$ and $\gamma_{i}^{-1}$ passing through
$y$ and $y'$, respectively. Clearly $h'=g(h^{-1})$. Now 
$s'=h\cdot\gamma'\cdot h'$ covers $s$ and so it is closed. Then we have
$h'=h^{-1}$ implying $g=1$ and $y=y'$. So we proved that $M_{1}$ is
trivial over each $\widetilde W_{i}$ and, in particular, it is trivial
over each $W_{i}$. This means that $M_{1}$ is defined by a locally constant
cocycle $\{c_{ij}\}$ (with values in $R_{1}$) defined on the cover 
$(W_{i})_{i\in I}$ of $M$. But then $\{c_{ij}\}$ is the pullback of a cocycle
defined on $(U_{i})_{i\in I}$ because the fibres of $f$ are connected.
This cocycle determines a principle bundle $C_{1}$ on $C$ with the fibre
$R_{1}$ such that $f^{*}C_{1}=M_{1}$. In particular, $C_{1}$ determines a
homomorphism $q:\pi_{1}(C)\rightarrow R_{1}$ such that $q\circ f_{*}$ is the
quotient homomorphism $\pi_{1}(M)\rightarrow R_{1}$. This implies that
$Ker(f_{*})\subset R$.\ \ \ \ \ $\Box$\\
{\bf 2.6.} We use use Lemma \ref{kern} to prove the following result.
\begin{Lm}\label{fini}
Let $\rho\in Hom(\pi_{1}(M),T_{2})$ be such that 
$Ker(f_{*})\subset Ker(\rho_{a})$ where $\rho_{a}$ is unitary.
Then $\rho(Ker(f_{*}))\subset N_{2}$ is a free abelian group of finite rank.
\end{Lm}
{\bf Proof.} Without loss of generality we will assume that 
$Ker(f_{*})\not\subset Ker(\rho)$.
Let $i:V\hookrightarrow\pi_{1}(M)$ be the fibre over a regular
value of $f$ and $H=i_{*}(\pi_{1}(V))\subset\pi_{1}(M)$ be the corresponding
normal subgroup. First, we prove that $\rho|_{H}$ is not trivial. Indeed,
let $E$ be the flat vector bundle associated to $\rho_{a}$. By our assumptions
$E=f^{*}L$ where $L$ is a flat vector bundle on $C$  with unitary structure
group. It is well known that the equivalence class of $\rho$ is defined by
a harmonic 1-form $\eta$ with values in $E$. Since the structure group of
$E$ is unitary, the Hodge decomposition implies that 
$\eta=\omega_{1}+\omega_{2}$, where $\omega_{1}$ is holomorphic and
$\omega_{2}$ is antiholomorphic. From the proof of Lemma \ref{van} we know 
that $\rho_{H}$ is defined by integration of $\eta$ along closed paths of 
$i(V)$. Assume to the contrary that $\rho|_{H}$ is trivial. This means that
$\eta|_{i(V)}$ represents 0 in $H^{1}(i(V),\Co)$ and so 
$\omega_{1}|_{i(V)}=\omega_{2}|_{i(V)}=0$. But then Lemma \ref{vanish}
implies that $\omega_{i}|_{V}=0$ $(i=1,2)$ for any fibre $V$ of $f$.
Let us consider a path $\gamma$ which represents an element of 
$\pi_{1}(V_{i},x_{i})$ (the notation is as in Lemma \ref{kern}). Let
$\gamma'=\gamma_{i}\cdot\gamma\cdot\gamma_{i}^{-1}$. Then the monodromy
argument shows that $\rho(\gamma')=C\cdot A\cdot C^{-1}$ where $C\in T_{2}$,
$A\in N_{2}$ and $A$ is defined by integration of $\eta$ along $\gamma$.
Therefore $\rho(\gamma')=1\in N_{2}$. Let $[\gamma']\in\pi_{1}(M)$ be the
element defined by $\gamma'$. According to Lemma \ref{kern},
$Ker(f_{*})$ is the minimal normal subgroup containing all such $[\gamma']$
(for any $i$). The above argument then shows that $\rho|_{Ker(f_{*})}$ is
trivial which contradicts to our assumption.

Next, we will prove that $\rho_{a}(\pi_{1}(M))$ consists of algebraic
integers. Let $H_{f}\cong\Z^{k}$ be the free part of $H^{ab}=H/DH$. In
what follows we identify $H_{f}$ with $\Z^{k}$. Let
$s:\pi_{1}(M)\rightarrow SL_{k}(\Z)$ be the representation induced by the
action of $\pi_{1}(M)$ on $H$ by conjugation.
Since $\rho(H_{f})=\rho(H)$ is non-trivial and consists of unipotent matrices,
$k\geq 1$ and identity (\ref{ident}) implies that the action of
$\rho(\pi_{1}(M))$ on $\rho(\Z^{k})$ by conjugation is defined 
by multiplication of non-diagonal elements of $\rho(\Z^{k})$ by elements of
$\rho_{a}(\pi_{1}(M))$. Let $e_{1},\dots,e_{k}$ be the standard basis
in $\Z^{k}$. Identifying $N_{2}$ with $\Co$
we obtain a nonzero vector $v=(\rho(e_{1}),\dots,\rho(e_{k}))\in\Co^{k}$.
Then a simple calculation shows that $v$ is a nonzero eigen vector of $s(g)$
with the eigen value $\rho_{a}(g)$, $g\in\pi_{1}(M)$. This completes the
proof.

Now let us finish the proof of the lemma. Let $\{g_{j}\}$, $1\leq j\leq l$,
be the {\em finite family} which consits of all generators of all $G_{i}$
(as in Lemma \ref{kern}). Then according to Lemma \ref{kern},
$\rho(Ker(f_{*}))\subset\rho(\pi_{1}(M))$ is the minimal normal subgroup
containing all $h_{j}:=\rho(g_{j})$. Since $N_{2}$ is abelian,
identity (\ref{ident}) shows that $\rho(Ker(f_{*}))$ is a discrete 
abelian group generated by all elements $\rho_{a}(g)\cdot g_{j}$,\
$g\in H_{1}(M,\Z)$, $1\leq j\leq l$.
Let $t_{1},\dots,t_{p}$ be the family of generators of $H_{1}(M,\Z)$.
Then the above arguments imply that each $\rho_{a}(t_{i}), 
\rho_{a}(-t_{i})$ 
$(i=1,\dots,p)$ are roots of polynomials of degree $k$ with integer 
coefficients and leading coefficients $1$. Since any $g\in H_{1}(M,\Z)$
can be written as 
$\sum_{i=1}^{p}a_{i}t_{i}+\sum_{i=1}^{p}b_{i}(-t_{i})$ with 
$a_{i},b_{i}\in\Z_{+}$, the previous shows that the abelian group 
$\rho(Ker(f_{*}))$ is generated (over $\Z$) by the elements 
$$
\rho_{a}(t_{1})^{a_{1}}\dots\rho_{a}(t_{p})^{a_{p}}\cdot
\rho_{a}(-t_{1})^{b_{1}}\dots\rho_{a}(-t_{p})^{b_{p}}\cdot g_{j}
$$
where $a_{i}, b_{i}$ are nonnegative integers satisfying
$0\leq a_{i}\leq k-1$, $0\leq b_{i}\leq k-1$, and $1\leq j\leq l$. Since
the number of such generators is finite, $\rho(Ker(f_{*}))$ is finitely
generated. It is also free as a subgroup of the free group $N_{2}$.

The proof of the lemma is complete.\ \ \ \ \ \ $\Box$\\
{\bf 2.7. Proof of Theorem \ref{curv}.} 
Let $\rho\in Hom(\pi_{1}(M),T_{2})$ be such that 
$Ker(f_{*})\subset Ker(\rho_{a})$ but $Ker(f_{*})\not\subset Ker(\rho)$.
According to Proposition \ref{normsubgr}, $\rho_{a}$ is unitary. Then
according to Lemma \ref{fini}, $\rho(Ker(f_{*}))\subset N_{2}$ is 
isomorphic to $\Z^{r}$, $r\geq 1$. We set $G=Ker(f_{*})\cap Ker(\rho)$ and
$R:=\pi_{1}(M)/G$. Then $R$ is defined by the sequence
$$
\{1\}\longrightarrow Ker(f_{*})/G\longrightarrow
R\longrightarrow\pi_{1}(C)\longrightarrow\{1\}\ . 
$$ 
where $Ker(f_{*})/G\cong\rho(Ker(f_{*}))=\Z^{r}$. Since 
$\rho(Ker(f_{*}))\subset N_{2}$, the above definition implies that
the natural action of $D\pi_{1}(C)$ on $Ker(f_{*})/G$ is trivial.
Thus according to Proposition \ref{curve}, all characters of the induced
representation $\pi_{1}(C)\rightarrow GL_{r}(\Co)$ are torsion. 
In particular, $\rho_{a}$ is torsion as one of such characters.

The proof of the theorem is complete.\ \ \ \ \ \ $\Box$

\end{document}